\documentclass{ifacconf}

\usepackage{amssymb}
\usepackage{amsmath}
\usepackage{graphicx}      % include this line if your document contains figures
\usepackage{natbib}        % required for bibliography
\usepackage{color}

\newcommand{\vect}[1]{\stackrel{\rightarrow}{#1}}

\newtheorem{theorem}{Theorem}
\newtheorem{remark}{Remark}
\newtheorem{lemma}{Lemma}

\newcommand{\be}{\begin{equation}}
\newcommand{\ee}{\end{equation}}
\newcommand{\ba}{\begin{array}}
\newcommand{\ea}{\end{array}}

\newcommand{\la}{\left\langle}
\newcommand{\ra}{\right\rangle}

\newcommand{\lb}{\lambda}

\newcommand{\al}{\alpha}

\newcommand{\tb}{{\overline{\tau}}}

\newcommand{\RR}{{\mathbf{R}}}

\begin{document}
\begin{frontmatter}

\title{Bounds on time-optimal concatenations of arcs for 
two-input driftless 3D systems 
\thanksref{footnoteinfo}} 

\thanks[footnoteinfo]{ The project leading to this publication has received funding from the ANR project SRGI ANR-15-CE40-0018.}

%N° 765267

\author[First]{M. Sigalotti}

\address[First]{Inria \& Sorbonne Universit\'e, Universit\'e de Paris, CNRS, Laboratoire Jacques-Louis Lions, Paris, France (e-mail: mario.sigalotti@inria.fr)}

\begin{abstract}                % Abstract of not more than 250 words.
We study a driftless system on a three-dimensional manifold driven by two scalar controls.
We assume that each scalar control has an independent bound on its modulus
and we prove that, locally around every point where the controlled vector fields 
satisfy some suitable nondegeneracy Lie bracket condition,
every time-optimal trajectory has at most 
five bang or singular arcs. The result is obtained using 
first- and second-order necessary conditions for optimality.  
\end{abstract}


\end{frontmatter}
%===============================================================================

\section{Introduction}

The regularity of optimal trajectories is a much studied problem in control theory, because 
the presence of discontinuities 
has obvious implications on the implementability and robustness of optimal feedbacks. 
The regularity of the value function is also strongly related with the regularity of the optimal trajectories (see, for instance, \cite{Schattlerbook}). When the value function encodes a distance in a length-space, the regularity of geodesics carries information on the  properties of 
balls and other relevant geometric objects \cite{agrachev_barilari_boscain_2019}.  

We study here a class of control systems which can be written in the form
\be\label{eq:fins0}
\dot q(t)=u_1(t) X_1(q(t))+u_2(t)X_2(q(t)),
\ee
where the state $q$ evolves on a smooth manifold $M$, $X_1$, $X_2$ are smooth vector fields, and the control $u=(u_1,u_2)$ takes values in the square $[-1,1]^2$. 
This setup corresponds to driftless two-input systems for which the modulus of each control parameter has an independent bound.
Under some natural Lie algebra rank condition, the time-optimal value function for \eqref{eq:fins0} can also be seen as the length distance associated with a sub-Finsler structure (see \cite{BBLDS}).

Few is known on the general  structure of time-optimal trajectories for systems of the type \eqref{eq:fins0} apart from some homogeneous cases studied in \cite{BreuillardLeDonne} (where $M$ is the Heisenberg group and $X_1$, $X_2$ are two left-invariant generators of the horizontal distribution), \cite{BBLDS} (where sub-Finsler versions  of the Grushin plane and the Martinet distribution have been considered), and \cite{Ardentovetal} (where the Cartan group is considered). 
The inhomogeneous 2D case has been studied in \cite{AliCharlot}.

More is known in the case where the control $u=(u_1,u_2)$ takes values in a ball, 
since this leads to the sub-Riemannian framework (see \cite{agrachev_barilari_boscain_2019,jean-book,rifford-book}). 
It should however be mentioned that even in this case the minimal regularity of time-optimal trajectories is still an open problem (see \cite{SSS,Hakavuori,Monti} for recent results and a state of the art about this longstanding question). 
When the control takes values in a ball a drift is added to the dynamics fewer results are known (see \cite{AgrachevBiolo1,AgrachevBiolo2,Caillauetal}). 

Another case where more results are available is when one of the two control inputs $u_1$ or $u_2$ is constant. In this case we recover the case of  a single-input control-affine system with control in a compact interval.
For such kind of systems on a two-dimensional manifold, the situation has been deeply analyzed and general results covering the generic and the analytic case can be found in \cite{Sussmann2Dgeneric, Sussmann2Danalytic}. A monograph where 2D optimal syntheses are studied in details is \cite{BoscainPiccoli}. 
The results in the 3D case are less complete and do not cover the generic case, but only the local behavior of time-optimal trajectories near points where the brackets of the vector fields $X_1$ and $X_2$ have degeneracies of corank 0 (i.e., near generic points) or 1 (see \cite{AG1,AgrachevSigalotti,KS89,schattler1988local,schattler88altro,sussmann1986envelopes}) and some but not all cases of corank larger than 1 (\cite{bressan1986,JDCS}). In \cite{JDCS},  the regularity of time-optimal trajectories near generic points in dimension 4 is also characterized.  

The main type of result contained in the papers mentioned above is a guarantee that small-time optimal trajectory are the concatenation of at most a given number of bang and singular arcs (with limitations on the possible concatenation). One of the relevant consequences of this kind of result is that they allow to rule out the appearance of the Fuller phenomenon (also called \emph{chattering}). This kind of rather radical singularity of optimal trajectory is known to be 
typical (i.e., not removable by small perturbations) in large enough dimension (see \cite{kupka1990ubiquity,ZelikinBorisov,Zelikin-multi}). Some restriction on the wildness of the Fuller behavior (in particular, on the iteration order of accumulations of switching times) can be found in \cite{BoarottoSigalotti,BoarottoChitourSigalotti}.

In this paper we provide a bound on the worst-regularity 
behavior of time-optimal trajectories of 
\eqref{eq:fins0} near generic points in the case where $M$ has dimension 3. Our main result, Theorem~\ref{t:main}, states that such trajectories are necessarily concatenations of bang and singular arcs, and that the number of such arcs is not larger than 5. It is interesting to observe that in the corresponding homogenous 3D case, the sub-Finsler structure on the Heisenberg group, small-time optimal trajectories with 5 arcs do exist (\cite{BreuillardLeDonne,BBLDS}). This means that  the bound that we provide in this paper is sharp. 
The main technical step in the proof of our main result used the second-order necessary conditions for optimality proposed in \cite{AG1}.

\section{Statement of the time-optimal problem and necessary conditions for optimality
}

Throughout the paper $M$ denotes a smooth (i.e., $\mathcal{C}^\infty$) complete manifold and 
$X_1,X_2$  are 
 two smooth vector fields on $M$. 
We associate with $M$, $X_1$, and $X_2$ the dynamics 
\be\label{fins}
\dot q= u_1 X_1(q)+u_2 X_2(q),\quad q\in M,\quad u_1,u_2\in[-1,1].
\ee
We also introduce the notations $u=(u_1,u_2)$ and  
\begin{align*}
X(u)&=u_1 X_1+u_2 X_2,\\
X_+&=X_1+ X_2,\\
X_-&=X_1-X_2.
\end{align*}

It is well known that the time-optimal trajectories of a control system of the form \eqref{fins} heavily depend on the commutativity properties of the vector fields $X_1$, $X_2$, which are infinitesimally characterized by the iterated Lie brackets  between them. 

Given two smooth vector fields $X$ and $Y$ we write $[X,Y]$ 
to denote the Lie bracket between $X$ and $Y$, and ${\rm ad}_X$ for the operator $[X,\cdot]$ acting on the space of smooth vector fields by
\[{\rm ad}_X(Y) =[X,Y].\]

To reduce the notational burden we also set 
\[X_{12}=[X_1,X_2],\quad X_{+12}=[X_+,X_{12}],\quad X_{-12}=[X_-,X_{12}].\]
The Lie algebra generated by $\{X_1,X_2\}$, denoted by ${\rm Lie}\{X_1,X_2\}$, is the minimal linear subspace of the space of all smooth vector fields on $M$ containing $\{X_1,X_2\}$ and invariant both for ${\rm ad}_{X_1}$ and 
${\rm ad}_{X_2}$.

The time-optimal control problem for \eqref{fins} consists in finding the trajectories of 
\eqref{fins} 
connecting two given points $q_0$ and $q_1$ of $M$ in minimal time. The existence of at least one trajectory connecting any pair of points and minimizing the time is a consequence of Chow's theorem and Filippov's theorem, under the assumption that 
\[\{V(q)\mid V\in {\rm Lie}\{X_1,X_2\}\}=T_q M\]
for every $q\in M$ (see, for instance, \cite{AgrachevSachkov}).

\subsection{First-order necessary  conditions for optimality: the Pontryagin maximum principle}

A well-known optimality condition satisfied by all time-optimal trajectories is the Pontryagin maximum principle (PMP, for short). In order to fix some notations, let us recall here its statement. 

Let $\pi:T^*M\to M$ be the cotangent bundle, and $s\in \Lambda^1(T^*M)$ be the tautological Liouville one-form on $T^*M$. The non-degenerate skew-symmetric form $\sigma=ds\in \Lambda^2(T^*M)$ endows $T^*M$ with a canonical symplectic structure.
With any smooth function $p:T^*M\to \RR$ let us associate its smooth Hamiltonian lift $\vect{p}\in
	{\cal C}(T^*M,TT^*M)$ by the condition
	\be\label{eq:Hamlift}
		\sigma_\lambda(\cdot,\vect{p})=d_\lambda p.
	\ee
	
Introducing the control-dependent Hamiltonian function ${\cal H}:\RR^{2}\times T^*M\to\RR$ by
	\be\label{eq:maxH}
		{\cal H}(\lambda,v)=\langle \lambda, v_1X_1(q)+v_2X_2(q) \rangle, \ \ q=\pi(\lambda),
	\ee
	the statement of the PMP for the time-optimal problem associated with system \eqref{fins} is the following (see, for instance, \cite{AgrachevSachkov}).
	
	\begin{thm}[PMP]\label{thm:PMP}
		Let $q:[0,T]\to M$ be a time-optimal trajectory of \eqref{fins}, associated with a 
		control $u(\cdot)$.
		Then there exists an absolutely continuous curve
		$\lambda:[0,T]\to T^*M$ such that
		$(q(\cdot),\lambda(\cdot),u(\cdot))$ is an \emph{extremal triple}, i.e., in terms of the control-dependent Hamiltonian $\mathcal{H}$ introduced in \eqref{eq:maxH}, one has 
		\begin{align}
			&\lambda(t)\in T^*_{q(t)}M\setminus \{0\},\quad \forall t\in [0,T],\\
			\label{eq:maxPMP}
			&		\mathcal{H}(\lambda(t),u(t))=\max\{\mathcal{H}(\lambda(t),v)\mid v\in[-1,1]^2\}\\
			&\quad  \ \mathrm{for\ a.e.\ }t\in[0,T],\nonumber\\
			&\dot{\lambda}(t)=\vect{\mathcal{H}}
			(\lambda(t),u(t)),\quad \mathrm{for\ a.e.\ }t\in[0,T].\label{eq:PMP}
		\end{align}
	\end{thm}
	
		For any extremal triple $(q(\cdot),\lambda(\cdot),u(\cdot))$, we call the corresponding trajectory $t\mapsto q(t)$ an \emph{extremal trajectory}, and the curve $t\mapsto \lambda(t)$ its associated \emph{extremal lift}.
	
	Let us introduce the smooth functions $\phi_1,\phi_2:T^*M\to\RR$ by
	\[
		\phi_i(\lambda)=\langle \lambda, X_i(q)\rangle,\quad q=\pi(\lambda), \quad i=1,2.
	\]
	
A maximal open interval of $[0,T]$ on which $\phi_1$ and $\phi_2$ are both different from zero is said to be a \emph{bang arc}. A  \emph{$u_1$-singular} (respectively, \emph{$u_2$-singular}) arc is a 
maximal open interval of $[0,T]$ on which $\phi_1\equiv 0$ while $\phi_2$ is different from zero (respectively, $\phi_2\equiv 0$ while $\phi_1$ is different from zero). 
An arc which is either $u_1$-singular or $u_2$-singular is said to be \emph{singular}. 
A point separating two bang arcs is said to be a \emph{switching time}.
We say that the control $u_i$ \emph{switches} at the switching time $\tau$ if $u_i$ (which is constant and of modulus $1$ on sufficiently small left- and right-neighborhoods of $\tau$) changes sign at $\tau$. In particular, if $u_i$ switches at $\tau$ then it follows from the maximality condition of the $PMP$ that $\phi_i(\tau)=0$.
A trajectory is said to be \emph{bang-bang} if it is the concatenation of finitely many bang arcs. 

A useful consequence of the PMP 
 is that for every 
smooth vector field $Y$ on $M$, for every extremal triple associated with \eqref{fins}
	the identity
	\be\label{eq:diffvectfield}
		\frac{d}{dt}\langle \lambda(t),Y(q(t)) \rangle=\langle \lambda(t),[X(u(t)),Y ](q(t))\rangle
	\ee
	holds true for a.e. $t$.

\subsection{Second-order optimality conditions}
We recall in this section a useful second-order necessary 
condition for a  trajectory with piecewise constant control to be 
time-optimal,  
obtained in \cite{AG1}
using time-reparameterizations as variations of the reference control signal. 

In order to state the theorem, let us introduce the following notation. Given a control $u:[0,T]\to [-1,1]^2$ and two times $s,t\in [0,T]$, denote by $P^{u}(s,t)$ the flow of \eqref{fins} from time 
$s$ to time $t$. Since we are interested in local properties, we can always assume that 
solutions of \eqref{fins} exist globally, which ensures that  $P^{u}(s,t)$ is a defined on the entire manifold $M$. As a consequence, by standard properties of solutions of ODEs, 
$P^{u}(s,t):M\to M$ is a diffeomorphism.

For every diffeomorphism $P:M\to M$, every point $q\in M$, and every $v\in T_qM$, $P_{*,q}v\in T_{P(q)}M$ denotes the push-forward of $v$ obtained by applying the differential of $P$ at $q$.  
For every diffeomorphism $P:M\to M$ and every vector field $Y$, the \emph{push-forward of $Y$ by $P$} is the vector field
\[P_*Y:q\mapsto P_{*,P^{-1}(q)}Y(P^{-1}(q)).\]

Using these notations, 
we can deduce from the PMP that, if $(q(\cdot),\lambda(\cdot),u(\cdot))$ is an extremal triple on $[0,T]$, then, for $i=1,2$ and for $t,\tb\in [0,T]$,
\be\label{eq:pmppush}
\phi_i(t)=\langle \lambda(\tb),(P^u(\tb,t)_* X_i)(q(\tb))\rangle.
\ee

\begin{theorem}\label{thm2nd}
Let  $q:[0,T]\to M$ be a time-optimal trajectory for \eqref{fins}
 and let
$u(\cdot)=(u_1(\cdot),u_2(\cdot))$ be the corresponding control function.  Assume that $u(\cdot)$ is piecewise constant on $[0,T]$, with $K$ non-removable discontinuities $\tau_1<\tau_2<\cdots<\tau_K$ in $(0,T)$.
 Denote by $u_0,\dots,u_k$ the successive values of $u(\cdot)$ on the $K+1$ bang arcs. 
Assume that $q(\cdot)$ admits a unique extremal lift $\lb(\cdot)$ up to multiplication by a positive scalar.
Fix $\tb$ in $[0,T]$ and let 
\be\label{h_i}
h_i=
P^u(\tb, \tau_i)_*X(u_i),
\ee
for $i=0,\dots,K$.
Let $Q$ be the quadratic form
\be\label{2nd} Q(\alpha)=\sum_{0\le i < j\le K}\alpha_i\alpha_j\la \lb(\tb),
[h_i,h_j](q(\tb))\ra
\,,
\ee
defined on the space
\be\label{space_of_alphas}
\begin{split}H=&
\left\{\alpha=(\alpha_0,\dots,\alpha_K)\in \RR^{K+1}\mid \right.\\
&\left.\sum_{i=0}^K\alpha_i=0,\ \sum_{i=0}^K \alpha_i h_i(q(\tb))=0\right\}.
\end{split}
\ee
Then $Q\leq 0$.
\end{theorem}

\begin{remark}
The theorem above can be extended from time-minimal to (locally) time-maximal trajectories.
This can be done following the lines of \cite[Theorem 2]{AgrachevSigalotti} where the notion of
quasi-optimal trajectory is introduced to cover both kind of properties. 
\end{remark}

\section{Bound on the number of arcs}

Let us assume from now on that $M$ has dimension $3$. 
This section contain the main results on the admissible concatenation of arcs, locally near a point where suitable bracket independence conditions are satisfied. 
We collect the main results in a single statement, Theorem~\ref{t:main} below.

\begin{theorem}\label{t:main}
Let $q_0\in M$ 
and $\Omega$ be a neighborhood of $q_0$, compactly contained in $M$,
such that 
$(X_1,X_2,X_{12})$, $(X_1,X_{12},X_{+12})$, $(X_1,X_{12},X_{-12})$, $(X_2,X_{12},X_{+12})$, and $(X_2,X_{12},X_{-12})$ are moving bases on the closure $\overline{\Omega}$ of $\Omega$. 
Then there exists $T>0$ such that every 
time-optimal
trajectory 
$q:[0,T']\to M$ of \eqref{fins} contained in $\Omega$ and such that $T'\le T$ 
is the concatenation of at most 5  bang or singular arcs. 
Moreover, if $q(\cdot)$ contains a singular arc, then it is the concatenation of at most a bang, a singular, and a bang arc.
\end{theorem}

The main technical step in the proof of the theorem is contained in the next lemma, which focuses on the situation in which both $u_1$ and $u_2$ switch along the time-optimal  trajectory. 
Notice that in this case 
we can relax the assumptions on the triples of vector fields which should be linearly independent on the considered neighborhood $\Omega$. 
Actually, we just need to assume that $X_1$, $X_2$, and $X_{12}$ are linearly independent on $\overline{\Omega}$. 

\begin{lemma}\label{l:bang}
Let $q_0\in M$ 
and $\Omega$ be a neighborhood of $q_0$, compactly contained in $M$,
such that 
$X_1$, $X_2$, and $X_{12}$ are linearly independent on the closure $\overline{\Omega}$ of $\Omega$. 
Then there exists $T>0$ such that every 
bang-bang time-optimal
trajectory 
$q:[0,T']\to M$ of \eqref{fins} contained in $\Omega$, undergoing switchings both in $u_1$ and in $u_2$, and such that $T'\le T$
is the concatenation of at most 5 
arcs. 
\end{lemma}
\begin{pf}
Consider a bang-bang extremal trajectory $q:[0,T']\to M$ of \eqref{fins} contained in $\Omega$.
Assume that $q(\cdot)$ is the concatenation of 6 bang arcs and that both $u_1$ and $u_2$ switch. 
We are going to show that, for $T'$ small enough, $q(\cdot)$ is not optimal.

Let $\lambda(\cdot)$ be an extremal lift of $q(\cdot)$ and define
\be\label{eq:phiss}
\phi_\star(t)=\la \lb(t),X_\star(q(t))\ra,\qquad t\in [0,T'],\quad \star\in\{1,2,12
\}.
\ee
It follows from \eqref{eq:diffvectfield} that
\begin{equation}\label{eq:derbra}
\dot \phi_\star(t) =\la \lb(t),[u_1(t)X_1+u_2(t) X_2,X_\star](q(t))\ra
\end{equation}
for $\star\in\{1,2,12\}$ and 
for almost every $t\in [0,T']$. In particular, 
\be\label{eq:derbraplus}
\begin{split}
\dot \phi_1(t) &=-u_2(t)\phi_{12}(t),\\ 
\dot \phi_2(t) &=u_1(t)\phi_{12}(t),
\end{split}
\ee
for almost every $t\in [0,T']$.

Assume for now that $\lb(\cdot)$ is normalized in such a way that 
\begin{equation}\label{eq:normalization}
\max\{|\phi_1(0)|,|\phi_2(0)|,|\phi_{12}(0)|\}=1.
\end{equation}

We first claim that, up to taking $T$ small enough, $\phi_{12}$ does not change sign on $[0,T']$. 
In order to prove the claim, notice that, since $\lb(\cdot)$ solves the time-dependent Hamiltonian system
\eqref{eq:PMP}
described by the Pontryagin maximum principle, its growth admits a uniform bound among extremal trajectories in $T^*\Omega$. 
It then follows from \eqref{eq:derbra} that there exist $T,C>0$ (independent of $q(\cdot)$) such that $|\dot \phi_\star(t)|\le C$
for $\star\in\{1,2,12\}$ and $t\in [0,T']$, provided that $T'<T$.
Because of the assumption that both $u_1$ and $u_2$ switch along the trajectory $q(\cdot)$, we know that both 
 $\phi_1$ and $\phi_2$ have a zero on $[0,T']$. 
This implies, in particular, that $|\phi_1|,|\phi_2|\le C T'$ on $[0,T']$. 
Let $T'$ be smaller that $1/C$.
It follows that 
$|\phi_1|,|\phi_2|<1$ on $[0,T']$, and hence,  from \eqref{eq:normalization}, 
one has 
$|\phi_{12}(0)|=1$. 
As a consequence, up to modifying $C$ uniformly with respect to $q(\cdot)$, 
$|\phi_{12}-1|\le CT'$ on $[0,T']$, and, in particular, $\phi_{12}$ has constant sign. This concludes the proof of the claim that $\phi_{12}$ can be assumed not to change sign on $[0,T']$.

Let us  focus now on the admissible concatenations of switches.
Assume for a moment that $u_2$ switches at the same time $\tau\in (0,T')$ at which 
$u_1$ does. Then $\dot \phi_1$ changes sign at $\tau$. This means that $\phi_1$ does not change sign at $\tau$, contradicting the assumption that $u_1$  switches at $\tau$. 
We can then assume that the switches of $u_1$ and $u_2$ are distinct.

Denote by $\tau_0\in (0,T')$
the first switching time  for $q(\cdot)$. 
Up to exchanging the roles of $X_1$ and $X_2$, we can assume that $\tau_0$ is a switching time for $u_1$. 
Up to exchanging $X_2$ and $-X_2$, we can assume that $u_2=-1$ in a neighborhood of $\tau_0$.
Denote by  $\tau_1$ the smaller switching time  for $q(\cdot)$ with $\tau_1>\tau_0$.
Assume for now that the (constant) sign of 
$\phi_{12}$ on $[0,T']$ is $-1$.
Hence,  
$\phi_2$ is increasing and $\phi_1$ decreasing on $(\tau_0,\tau_1)$. This implies that $\tau_1$ is a switching time for $u_2$. The same reasoning allows to continue the argument and deduce that the sequence of constant controls $(u_1,u_2)$ for $q(\cdot)$ on its 6 bang arcs follows 
the (periodic) 
pattern 
\begin{equation}\label{pattern}
\begin{pmatrix}1\\-1\end{pmatrix}\stackrel{\to}{\tau_0} \begin{pmatrix}-1\\-1\end{pmatrix}\stackrel{\to}{\tau_1} \begin{pmatrix}-1\\1\end{pmatrix}\stackrel{\to}{\tau_2} \begin{pmatrix}1\\1\end{pmatrix}\stackrel{\to}{\tau_3}\begin{pmatrix}1\\-1\end{pmatrix} \stackrel{\to}{\tau_4} \begin{pmatrix}-1\\-1\end{pmatrix},
\end{equation}
where we have denoted by $\tau_0,\dots,\tau_4$ the 5 switching times of $q(\cdot)$. The 6 bang arcs correspond then to the intervals $I_1=(0,\tau_0)$, $I_2=(\tau_0,\tau_1)$, and so on, up to $I_6=(\tau_4,T')$.

In the case  where $\phi_{12}>0$ on $[0,T']$, one analogously shows that the pattern is
\[\begin{pmatrix}-1\\-1\end{pmatrix}\stackrel{\to}{\tau_0}\begin{pmatrix}1\\-1\end{pmatrix} \stackrel{\to}{\tau_1} \begin{pmatrix}1\\1\end{pmatrix} \stackrel{\to}{\tau_2} \begin{pmatrix}-1\\1\end{pmatrix}\stackrel{\to}{\tau_3} \begin{pmatrix}-1\\-1\end{pmatrix}\stackrel{\to}{\tau_4} \begin{pmatrix}1\\-1\end{pmatrix}
\]
and the argument below can be adapted easily.

Let us denote by  $t_i=\tau_{i+1}-\tau_i$, $i=0,1,2,3$, the length of the four bang arcs connecting two switching times.

We are going to apply Theorem~\ref{thm2nd} to the trajectory $q(\cdot)$ taking $\tb=\tau_2$. 
Notice that $\phi_1(\tb)=0$, since $u_1$ is switching at $\tb$.

It is convenient to consider a new normalization of the covector $\lambda(\cdot)$ by imposing that 
\be\label{eq:normla}
\phi_{12}(\tb)=-1.
\ee 

One of the assumptions of Theorem~\ref{thm2nd} is the uniqueness of the extremal lift for $q(\cdot)$. In order to prove uniqueness, notice that $\lambda(\cdot)$ depends only on the value $\lambda(\tb)$ and that 
\[\langle \lambda(\tb),X_1(q(\tb))\rangle=0,\]
because $\tb$ is a switching time for $u_1$, and 
\[\langle \lambda(\tb),X_{12}(q(\tb))\rangle=-1,\]
because of the normalization of $\lb(\cdot)$. 
We are left to prove the uniqueness of the component $\phi_2(\tb)$ of 
 $\lambda(\tb)$ along  $X_{2}(q(\tb))$.

Using \eqref{eq:pmppush}, let 
us consider, for $t\in (0,t_2)$, the development 
\begin{align}
\phi_2(\tb+t)&=\langle \lambda(\tb),(
{e^{t(X_1+X_2)}}_*X_2
)(q(\tb))\rangle\nonumber\\
&=
\phi_2(\tb)a(q(\cdot),t)+b(q(\cdot),t),\label{eq:ab}
\end{align} 
where the functions $a$ and $b$ depend on $t$ and on the 
bang-bang trajectory $q(\cdot)$,
but not on the extremal lift $\lambda(\cdot)$. 

We write below  
$O(t)$ to denote a quantity which can be bounded from above by a term of the form $C|t|$, 
for $t$ small enough, with $C$ uniform with respect to $q(\cdot)$.

With this notation, 
\[a(q(\cdot),t)=1+O(t),\quad b(q(\cdot),t)=O(t).\]
Since $\phi_2(\tb+t_2)=\phi_2(\tau_3)=0$, we deduce from \eqref{eq:ab} that $\phi_2(\tb)$ is a function of $q(\cdot)$ only (including the dependence on $t_2$), i.e., that $\phi_2(\tb)$ is uniquely identified by $q(\cdot)$. This completes the proof of the uniqueness of the extremal lift.

Pushing the computations a step further and noticing that 
\[
{e^{t_2(X_1+X_2)}}_*X_2
=X_2+t_2[X_1,X_2]+O(t_2^2)\]
on $\Omega$, 
we have that
\[a(q(\cdot),t_2)=1+O(t_2^2),\quad b(q(\cdot),t_2)=O(t_2^2).\]
Together with \eqref{eq:ab} evaluated at $t=t_2$, this yields 
\[ \phi_2(\tb)=t_2+O(t_2^2).\]

Repeating the argument on the interval $[\tau_1,\tb]$, 
we have that
\begin{align*}
\phi_2(\tb+t)&=\langle \lambda(\tb),
({e^{-t(-X_1+X_2)}}_*X_2)
(q(\tb))\rangle
\\
&=
\phi_2(\tb)(1+O(t^2))+O(t^2),
\end{align*} 
for $t\in [-t_1,0]$, yielding
\[ t_2+O(t_2^2)=\phi_2(\tb)=t_1+O(t_1^2).\]
Henceforth, 
\[ t_2=t_1+O(t_1^2).\]

Similarly, 
\[t_0,t_3=t_1+O(t_1^2),\]
and we then deduce that $T'=O(t_1)$. 
In particular, 
\[\max_{t\in [0,T'],\; i=1,2}|\phi_{i}(t)|=
O(t_1).\]

According to the relation $T'=O(t_1)$, we are left to prove that $q(\cdot)$ is not optimal if $t_1$ is small enough.

Using the notation from Theorem~\ref{thm2nd}, we have 
\begin{align*}
h_0&= {e^{-t_1(-X_-)}}_*{e^{-t_0 (-X_+)}}_*(X_-)=X_- +O(t_1),\\
h_1&= {e^{-t_1 (-X_-)}}_*(-X_+)=-X_++O(t_1),\\
h_2&=-X_-,\\
h_3&
=X_+,\\
h_4&={e^{t_2 (X_+)}}_*(X_-)=X_-+O(t_1),\\
h_5&={e^{t_2 (X_+)}}_*{e^{t_3 (X_-)}}_*(-X_+)=-X_++O(t_1).
\end{align*}
We can then evaluate
\[\sigma_{ij}=\la \lb(\tb),[h_i,h_j](q(\tb))\ra\]
for $0\le i<j\le 4$, obtaining 
\begin{align*}
\sigma_{01}&=2+O(t_1),
&\sigma_{02}&=O(t_1),\\
\sigma_{03}&=-2+O(t_1),
&\sigma_{04}&=O(t_1
),\\%\\
\sigma_{05}&=2+O(t_1),
&\sigma_{12}&=2+O(t_1),\\
\sigma_{13}&=O(t_1),
&\sigma_{14}&=-2+O(t_1),\\
\sigma_{15}&=O(t_1
),
&\sigma_{23}&=2,\\%\\
\sigma_{24}&=O(t_1),
&\sigma_{25}&=-2+O(t_1),\\
\sigma_{34}&=2+O(t_1),
&\sigma_{35}&=O(t_1),\\
\sigma_{45}&=2+O(t_1). 
\end{align*}

Decomposing the relation $\sum_{i=0}^5 \alpha_i h_i(q(\tb))=0$ on the basis $X_+(q(\tb)),X_-(q(\tb)),X_{12}(q(\tb))$ and collecting the components along $X_+(q(\tb))$ and $X_-(q(\tb))$, we get
\begin{align*}
0=&\alpha_0 O(t_1)-\alpha_1(1+O(t_1))+\alpha_3\\
&+\alpha_4 O(t_1)-\alpha_5(1+O(t_1)),\\
0=&\alpha_0(1+O(t_1))+\alpha_1 O(t_1)-\alpha_2\\
&+\alpha_4(1+O(t_1))+\alpha_5 O(t_1).
\end{align*}
Considering, in addition, the relation $\sum_{i=0}^5 \alpha_i=0$ and solving with respect to $\alpha_0,\alpha_1,\alpha_2$, we get 
\begin{align*}
\alpha_0&=-\alpha_3-\alpha_4+O(t_1;\alpha_3,\alpha_4,\al_5),\\%
\alpha_1&=\alpha_3-\al_5+O(t_1;\alpha_3,\alpha_4,\al_5),\\
\alpha_2&=-\alpha_3+O(t_1;\alpha_3,\alpha_4,\al_5),
\end{align*}
where $(\alpha_3,\alpha_4,\al_5)\mapsto O(t_1;\alpha_3,\alpha_4,\al_5)$ 
denotes a linear function whose coefficients are $O(t_1)$ in the sense introduced above.

We claim now that the space $H$ defined in \eqref{space_of_alphas}, taken here with $K=5$, is of dimension $3$ in $\RR^6$.
Indeed, if a vector $p$ is orthogonal to 
\[W=\left\{\sum_{i=0}^5 \alpha_i h_i(q(\tb))\mid \sum_{i=0}^5 \alpha_i=0\right\},\] 
then $p$ 
annihilates, in particular, $h_{i+1}(q(\tb))-h_i(q(\tb))$ for $i=0,\dots,4$. The same computations as those used to prove the uniqueness of the extremal lift of $q(\cdot)$ show that the orthogonal to 
$W$ is one-dimensional, that is, $W$ is 2-dimensional. Hence, the kernel $H$ of the map 
\[(\alpha_0,\dots,\alpha_5)\mapsto \sum_{i=0}^5 \alpha_i h_i(q(\tb))\]
on the 5-dimensional space 
\[\left\{(\alpha_0,\dots,\alpha_5)\mid \sum_{i=0}^5\alpha_i=0\right\}\]
 is of dimension 
$5-\dim W=3$.

The quadratic form $Q$ from Theorem~\ref{thm2nd} is then 
\[Q(\alpha_3,\alpha_4,\al_5)=
(\alpha_3,\alpha_4,\al_5) (M+O(t_1)) (\alpha_3,\alpha_4,\al_5)^T,\]
where 
\[ M=\begin{pmatrix} -4&0&2\\ 0&0&2\\2&2&0\end{pmatrix}.\]
Since $M$ has one positive and two negative eigenvalues, we conclude from Theorem~\ref{thm2nd} that the trajectory $q(\cdot)$ is not optimal for $t_1$ small enough.
This concludes the proof of Lemma~\ref{l:bang}.
\end{pf}

We are now ready to conclude the proof of Theorem~\ref{t:main}.

%\begin{pf}[Proof of Theorem~\ref{t:main}]
{\bf Proof of Theorem~\ref{t:main}.}
The proof works by considering separately several types of time-optimal trajectories $q:[0,T']\to \Omega$ of \eqref{fins}.

First consider the  case where 
either $\phi_1$ or $\phi_{2}$ never vanish on $[0,T']$. 
Up to exchanging the roles of $X_1$ and $X_2$, let us
assume that 
$\phi_2$ does not change sign and, consequently, 
$u_2$ is constantly equal to $+1$ or $-1$ on $[0,T']$. 

This means that the trajectory $q(\cdot)$ is time-optimal also for the single-input control-affine system
\be\label{eq:singinp}
\dot q=f(q)+u g(q),\quad u\in [-1,1],
\ee
with $g=X_1$ and $f$ equal either to $X_2$ or to $-X_2$. 
In particular $g,[f,g]$, and $[f+ g,[f,g]]$ are linearly independent on $\bar \Omega$ and the same is true for $g,[f,g]$, and $[f- g,[f,g]]$. 
We can deduce from \cite[Theorem 3]{AgrachevSigalotti} that 
$q(\cdot)$ 
is the concatenation of at most 3 bang arcs or a bang, a singular, and a bang arc.

Consider now the case where both 
$\phi_1$ and $\phi_{12}$ have at least one zero  on $[0,T']$. We claim that in this case 
$\phi_2$ is never zero on $[0,T']$, and the conclusion then follows by the case just considered.  
In order to prove the claim, normalize $\lb(\cdot)$ in such a way that 
\[\max(|\phi_1(0)|,|\phi_2(0)|,|\phi_{12}(0)|)=1.\]
Since 
the growth of $\lambda$ can be uniformly bounded on 
$\Omega$, we have that $|\phi_1|$ and $|\phi_{12}|$ can be bounded by $CT'$ for some positive constant $C$ independent on the trajectory $q(\cdot)$. 
By taking $T'$ small enough, 
it follows that 
$|\phi_2(0)|=1$.
Since, moreover, according to \eqref{eq:derbraplus},
\[|\dot \phi_2|\le |\phi_{12}|\le CT',\]
we can conclude that 
$\phi_2$ does never vanish, as claimed. 

The case where both 
$\phi_2$ and $\phi_{12}$ vanish on $[0,T']$ being completely analogous, we can assume from now on that
$\phi_{12}$ never vanishes on $[0,T']$.

We claim that in this case $q(\cdot)$ is bang-bang. 
Indeed, consider first the case where $q(\cdot)$ has a $u_1$-singular arc $(\tau_0,\tau_1)$. 
We deduce from the expression of $\dot\phi_1$ (cf. \eqref{eq:derbraplus}) and the nonvanishing of 
$\phi_{12}$ 
that $u_2$ must vanish on $(\tau_0,\tau_1)$. This, in turn, implies that also $\phi_2$ vanishes on $(\tau_0,\tau_1)$. 
Notice now that, by \eqref{eq:derbraplus}, $\dot \phi_2=u_1\phi_{12}$ which would also imply that $u_1\equiv 0$ on $(\tau_0,\tau_1)$. But, clearly, the control corresponding to a time-optimal trajectory cannot vanish on a nontrivial interval, invalidating the assumption that $q(\cdot)$ has a $u_1$-singular arc. 
The case of a $u_2$-singular arc is completely symmetric. 

In order to prove that $q(\cdot)$ is bang-bang we can then assume that there exists an open interval $(\tau_0,\tau_1)$ contained in $[0,T']$ on which $\phi_1$ and $\phi_2$ are different from zero. Assume that 
$(\tau_0,\tau_1)$ is maximal with this property and that $\tau_1<T'$. Then either $\phi_1$ or $\phi_2$ vanish at $\tau_1$. If they both vanish at $\tau_1$, then we deduce from \eqref{eq:derbraplus} that if one of the two functions $\phi_i$ changes sign, then also the derivative of the other function $\phi_{3-i}$ changes sign, which means that ${\rm sign}(\phi_{3-i})$ and hence $u_{3-i}$ stay constant. 
Moreover, since the sign of $\dot \phi_1$ and $\dot\phi_2$ do not change on a bang arc, the interval $(\tau_1,T')$ will be a bang arc. If they do not both vanish at the same point, then an analogous reasoning shows that $q|_{(\tau_1,T')}$ is either a single bang arc of that $q_{(\tau_0,T')}$ is the concatenation of bang arcs undergoing switches both in $u_1$ and $u_2$. 
Reasoning similarly for the interval $[0,\tau_0]$, we conclude the proof of the claim that $q(\cdot)$ is bang-bang. 

Moreover, the argument above shows that, under the condition that $\phi_12$ never vanishes,  the bang-bang trajectory $q(\cdot)$  is  either
the concatenation of at most 3 bang arcs, or it undergoes switches both in $u_1$ and $u_2$.
The  conclusion of the proof of Theorem~\ref{t:main} then follows from Lemma~\ref{l:bang}.

%\end{pf}

\section{Conclusion}

We considered in this paper the class of time-optimal driftless two-input control problems in  which both scalar control has an independent bound on its modulus. When the system is defined on a three-dimensional manifold and the controlled vector fields satisfy some generic Lie bracket independent condition at a given point, we prove that all small-time optimal trajectories near such a point are the concatenation of  at most five bang and singular arcs. 

The proof of this fact extensively uses the Hamiltonian formalism provided by the Pontryagin maximum principle and its related second-order sufficient conditions for  optimality. 

The proposed bound is sharp, as it has been illustrated by previous works in the literature considering the homogeneous case of the Heisenberg group, also known as Brockett integrator in the control literature.

\bibliography{biblio}

%% There are a number of predefined theorem-like environments in
%% ifacconf.cls:
%%
%% \begin{thm} ... \end{thm}            % Theorem
%% \begin{lem} ... \end{lem}            % Lemma
%% \begin{claim} ... \end{claim}        % Claim
%% \begin{conj} ... \end{conj}          % Conjecture
%% \begin{cor} ... \end{cor}            % Corollary
%% \begin{fact} ... \end{fact}          % Fact
%% \begin{hypo} ... \end{hypo}          % Hypothesis
%% \begin{prop} ... \end{prop}          % Proposition
%% \begin{crit} ... \end{crit}          % Criterion

\end{document}